\documentclass[12pt,reqno]{amsart} 

\usepackage[all]{xy}
\usepackage{supertabular}

\usepackage{bbm,amsmath}
\newcommand{\C}{\mathbbm{C}}
\DeclareMathOperator{\im}{im}

\newtheoremstyle{varthm-roman}% name
  {0.5ex}%      Space above, empty = `usual value'
  {0ex}%      Space below
  {\normalfont}% Body font
  {}%         Indent amount (empty = no indent, \parindent = para indent)
  {\bfseries}% Thm head font
  {:}%        Punctuation after thm head 
  {.5em}%     Space after thm head: " " = normal interword space;
  {\thmnote{#3}}% Thm head spec
\theoremstyle{varthm-roman}
\newtheorem*{varthm-roman}{}% all text supplied in the note
\newtheoremstyle{varthm-italic}% name
  {0.5ex}%      Space above, empty = `usual value'
  {0ex}%      Space below
  {\itshape}% Body font
  {}%         Indent amount (empty = no indent, \parindent = para indent)
  {\bfseries}% Thm head font
  {:}%        Punctuation after thm head 
  {.5em}%     Space after thm head: " " = normal interword space;
  {\thmnote{#3}}% Thm head spec
\theoremstyle{varthm-italic}
\newtheorem*{varthm-italic}{}% all text supplied in the note

\begin{document}

   \parindent0cm

   \title[Injective Analytic Maps]{Injective Analytic Maps - A
     Counterexample to the Proof}
   \author{Thomas Keilen}
   \address{Universit\"at Kaiserslautern\\
     Fachbereich Mathematik\\
     Erwin-Schr\"odinger-Stra\ss e\\
     D -- 67663 Kaiserslautern
     }
   \email{keilen@mathematik.uni-kl.de}
   \urladdr{http://www.mathematik.uni-kl.de/\textasciitilde keilen}
   \author{David Mond}
   \address{Mathematics Institute,
     University of Warwick,
     Coventry CV4 7AL,
     UK
     }
   \email{mond@maths.warwick.ac.uk}

   \date{September, 2003.}

     
   \begin{abstract}
     In \cite{Nem93} the author translates a conjecture of Le Dung
     Trang on the non-existence of injective analytic maps
     $f:\big(\C^n,0\big)\rightarrow\big(\C^{n+1},0\big)$ with $df(0)=0$
     into the non-existence of a hypersurface germ in
     $\big(\C^{n+1},0\big)$ with rather unexpected properties. However,
     the proof given in \cite{Nem93} contains an apparently fatal error,
     as we demonstrate with an example.
   \end{abstract}

   \maketitle

   In \cite{Nem93} the author addresses the problem whether the
   differential $df(0)$ of an
   injective analytic map germ $f:\big(\C^n,0\big)\rightarrow
   \big(\C^{n+1},0\big)$ can possibly be of rank less than $n-1$. 
   A long standing conjecture of Le Dung Trang for the case $n=2$ 
   states that this cannot be the case, even though it is
   not at all obvious how the topological fact of injectivity and the analytic
   datum on the rank of the derivative might relate to each other. 
   Analysing the image
   $(X,0)$ of $f$ as an analytic subspace of $\big(\C^{n+1},0\big)$,
   the author claims that a counter example to Le's conjecture would
   have an unexpected ``bad'' property.
   More precisely, he defines what it means for $(X,0)$ to be
   ``good'', and sets out to show that if $X$ is good then 
   the rank of $df(0)$ is
   at least $n-1$ and $(X,0)$ is an equisingular family of plane
   curves. However, the proof of this theorem contains a fundamental
   error, which -- as we are convinced after discussions with the author --
   cannot be repaired. We will outline the main ideas of the proof and
   give an example which shows that it does not work as described, and
   where it goes wrong. In order to keep the notation
   simple we restrict ourselves to the case where $n=2$. 

   We would like to point out that our example is not a 
   counter-example to the statement of the Theorem in \cite{Nem93} nor do we
   know of any such. It shows merely that the proof is wrong.

   Let us now recall the necessary definitions from \cite{Nem93}.

   \begin{varthm-roman}[Definition]
     A two-dimensional subgerm $(X,0)\subset\big(\C^3,0\big)$ is
     called \emph{good} if there exist coordinates $(w_1,w_2,w_3)$ for
     $\big(\C^3,0\big)$ and a map germ
     $F:\big(\C^3,0\big)\rightarrow(\C,0)$ defining $(X,0)$, i.~
     e.\ $X=F^{-1}(0)$, such that
     $W_0=X\cap \{w_1=0\}$ is an isolated plane curve
       singularity, and
     $\frac{\partial F}{\partial w_1}\not\in\left\langle w_1, 
       \frac{\partial F}{\partial w_2},\frac{\partial F}{\partial w_3}
     \right\rangle$.
   \end{varthm-roman}

   Nemethi then states the following

   \begin{varthm-italic}[``Theorem'']
     If the image $(X,0)$ of an injective analytic map germ 
     $f:\big(\C^2,0\big)\rightarrow \big(\C^3,0\big)$ is
     good, then the rank of $df(0)$ is at least one. Moreover, $(X,0)$ is an
     equisingular family of plane curve singularities over the base $(\C,0)$.
   \end{varthm-italic}
   
   The idea of the proof is to compare the two 
   isolated plane curve singularities
   $V_0=f_1^{-1}(0)$ and $W_0=X\cap \{w_1=0\}=\psi(V_0)$, where $f_i=w_i\circ
   f:\big(\C^2,0\big)\rightarrow (\C,0)$ for $i=1,2,3$ and
   $\psi=(f_2,f_3):\big(\C^2,0\big)\rightarrow \big(\C^2,0\big)$. 
   The Milnor fibre $V_t=f^{-1}(t)$ for $t\not=0$ maps via $\psi$ to
   $V_t'=\psi(V_t)$, which is in general singular. If $f$ is injective, then
   the restriction of $\psi$ to each level set of $f_1$ (i.e. to $V_t$)
   must also be injective.
   The vanishing cycles of $V_t$ must therefore be mapped homeomorphically
   by $\psi$ to non-trivial
   cycles in $V_t'$.
   Nemethi claims that under these circumstances, 
   the vanishing cycles of $V_t$, mapped by $\psi$ into $V_t'$, together 
   with the vanishing cycles of the singularities of $V_t'$ 
   (which it has acquired under the
   map $\psi$) together make up a complete set of vanishing cycles
   of a Milnor fibre of $W_0$. In $V_t$ one can choose vanishing cycles
   which do not pass through the (isolated) non-immersive points  
   of $\psi$.  In a smoothing of the singularities of $V_t'$, the 
   vanishing cycles can be confined to arbitrarily small neighbourhoods 
   (in the ambient space) of the points being smoothed, and thus 
   the vanishing cycles coming from the singularities of $V_t'$ have zero 
   intersection number with the images under $\psi$ of the vanishing 
   cycles coming from $V_t$. 
   {\it This implies that the Dynkin diagram of the
     isolated plane curve singularity $W_0$ is disconnected,} contradicting a
   well-known theorem of Lazzeri (\cite{Laz73}). 

   From this Nemethi concludes that
   one of the two sets of vanishing cycles must be empty, and thus that
   either $V_0$ or $V_t'$ is smooth. in the first case, the derivative at $(0,0)$
   of $f_1$ is not zero, and so the derivative of $f$ itself is not zero.
   In the second case, $V_t'$ is a Milnor fibre for $W_0$, and so $W_0$ and
   $V_0$ have the same Milnor number, from which it follows that $\psi$
   gives an isomorphism $V_0\to W_0$.
   From this Nemethi is able to 
   show that the germ $(X,0)$ is not good. 
   
   To make this argument rigorous, Nemethi has to show that the two types of
   cycles together really do form a basis of vanishing cycles in 
   a Milnor fibre of $W_0$.
   To do this he considers the deformation of $V_0$ induced by 
   $f_{1}:(\C^2,0)\to (\ell,0)=(\C,0)$.
   The image of this deformation under $\psi$ then gives a
   deformation of $W_0$ which can be induced from an ${\mathcal{R}}$-miniversal
   deformation $\Theta$ of $F_|:\{w_1=0\}\rightarrow(\C,0)$
   via base change $r$. The author claims then
   that a small perturbation 
   %%of $l$ and 
   of $r(\ell)$ gives rises to a Milnor fibre of $W_0$ in which the set of
   vanishing cycles splits into those
   coming from a Milnor fibre of $V_0$ and those arising from the
   singularities of $V_t'$.
   For this to be the case, it must be possible to deform $\ell'=r(\ell)$ 
   in a family to $\{\ell'_t\}_{t\in\C,0)}$ 
   in such a way that for $t\neq 0$, $\ell'_t$ intersects the discriminant
   $D$ in the base of the deformation $\Theta$ transversally in a finite 
   number of points, and that $\ell'_t\cap D$ does not
   meet the boundary of a good representative of the deformation. 
   The problem with the argument is that if $r(\ell)$ 
   is {\it contained} in $D$, then this is not
   in general possible. 
   And this is exactly what happens in our
   example, even though to see this one has to follow the
   constructions in the proof of the theorem very closely. For the
   details we refer to \cite{Kei93}. 

   An easy way
   to see that the proof must go wrong somewhere
   is to consider the following example. 
   \begin{displaymath}
     f:\big(\C^2,0\big)\rightarrow(\C^3,0\big):(x,y)\mapsto\big(y^3+x^2,x,y^2\big).
   \end{displaymath}
   Obviously $f$ is injective and
   \begin{displaymath}
     F:\big(\C^3,0\big)\rightarrow(\C,0):(w_1,w_2,w_3)\mapsto \big(w_1-w_2^2\big)^2-w_3^3
   \end{displaymath}
   is a defining equation of $(X,0)=\big(\im(f),0\big)$. In
   this case
   \begin{displaymath}
     V_0=f_1^{-1}(0)=\big\{y^3+x^2=0\big\}
   \end{displaymath}
   is a cusp, hence in particular not smooth, while
   \begin{displaymath}
     W_0=X\cap\{w_1=0\}=\big\{w_2^4-w_3^3=0\big\}
   \end{displaymath}
   is an $E_6$-singularity. Even though $f$ is injective, $V_0$ and $W_0$ 
   do not have the same Milnor number!

%   \bibliographystyle{amsalpha}
%   \bibliography{bibliographie}

\end{document}